\documentclass[secthm]{elsart}

 \usepackage{epsfig}
\usepackage{amssymb}
\usepackage{amsmath,amsfonts}

\def\C{{\cal C}}

\def\E{{\cal E}}

\begin{document}

\begin{frontmatter}



\title{Sofic groups and profinite topology on free groups \thanksref{thk}}
\thanks[thk]{The authors would like to thank
E. Gordon
for useful discussions. The work was partially supported by CONACyT grant
SEP-25750, and PROMEP grant UASLP-CA-21.}

\author{Lev Glebsky \corauthref{cor1}},
\ead{glebsky@cactus.iico.uaslp.mx}
\corauth[cor1]{Corresponding author}
\author{Luis Manuel Rivera}
\address{IICO-UASLP, Av. Karakorum 1470, Lomas 4ta Secci\'{o}n, San Luis Potos\'{\i},
SLP 7820 M\'{e}xico. Phone: 52-444-825-0892 (ext. 120)}

\begin{abstract}
We give a definition of weakly sofic groups (w-sofic groups). Our
definition is a rather natural extension of the definition of sofic
groups where instead of the Hamming metric  on symmetric groups we
use general bi-invariant metrics on finite groups. The existence of
non w-sofic groups is equivalent to a profinite topology property of
products of conjugacy classes in free groups.
\end{abstract}

\end{frontmatter}

\section{Introduction}
\label{intro}

The notion of sofic groups was introduced in \cite{grom,w} in
relation with the problem ``If for every group $G$ the injectivity
of cellular automata over $G$ implies their surjectivity?'' due to
Gottschalk \cite{Got}. The problem is still open and the class of
sofic groups is the largest class of groups for which the problem is
proved to have a positive solution \cite{w}, see also
\cite{CST_CM2007a,CST_CM2007b}. Sofic groups turned out to be
interesting from other points of view. For example, the "Connes
Embedding Conjecture" and "Determinant Conjecture" was proven for
sofic groups, \cite{es1}. The class of sofic groups is closed with
respect to various group-theoretic constructions: direct products,
some extensions, etc., see \cite{es1, es2}. It is still an open
question if there exists a non sofic group.

In the article we introduce the apparently more general class of
weakly sofic groups (w-sofic groups) (a sofic group is a w-sofic
group). We prove that the notion of w-sofic group is closely related
to some properties of the profinite topology on free groups. The
profinite topology on free groups (a base of neighborhoods of the
unit element consists of the normal subgroups of finite index)
 draw
considerable attention and one of the remarkable results here is
that a product of a finite number of finitely generated subgroups is
closed in the profinite topology (as well as in the pro-$p$
topology) \cite{Steinberg2005,BHDL,RZ}.

The question related to w-sofic groups is about the closure of
products of conjugacy classes in the profinite topology.

To finish this section, let us list some notation, that will be used
throughout the article. $S_n$ will denote the symmetric group on n
elements (the group of all permutations of a finite set
$[n]=\{1,...,n\}$). On $S_n$ we define the normalized Hamming metric
$h_n(f,g)=\frac{\mid \{a:(a)f \neq (a)g\} \mid}{n}.$ It is easy to
check that $h_n(\cdot,\cdot)$ is a  bi-invariant metric on $S_n$
\cite{MOP}. Sometimes we will omit the subscript $n$ in $h_n$. For a
group $G$ let $e_G$ denote the unit element of $G$. For
$A,B\subseteq G$ let, as usual, $AB=\{xy\;|\;x\in A,\;y\in B\}$. $X
< G$ and $X \vartriangleleft G$ will denote "$X$ is a subgroup of
$G$" and "$X$ is a normal subgroup of $G$", respectively. For $g\in
G$ let $[g]^G$ denote the conjugacy class of $g$ in $G$. For $g\in
G$ and $N<G$ let $[g]_N=gN$ denote the left $N$-coset of $g$.

\section{Profinite topology}

Closed sets in the profinite (the pro-$p$) topologies can be
characterized by the following separability property: a set
$X\subseteq G$ is closed in the profinite (the pro-$p$) topology iff
for any $g\not\in X$ there exists a homomorphism $\phi$ from $G$ to
a finite group (a finite $p$-group) such that
$\phi(g)\not\in\phi(X)$. It is clear, that for a closed $X\subset G$
and $g_1,g_2,...,g_k\not\in X$ there exists a homomorphism $\phi$
(to a finite group or to a finite $p$-group, depending on the
topology considered) such that
$\phi(g_1),\phi(g_2),...,\phi(g_k)\not\in \phi(X)$. (Just take
$\ker{\phi}=\cap\ker{\phi_i}$, where $\phi_i(g_i)\not\in\phi_i(X)$
). By the same way one can characterize the closure in the profinite
(the pro-$p$) topology. This characterization will be used in the
proof of Theorem~\ref{th_i_1}.

In the present article we are interested in the profinite topology
on a finitely generated free group $F$. It is known, that $[g]^F$ is
closed in the  pro-$p$ (and the profinite) topology for any $g\in
F$, (a conjugacy class in a free group is separable by homomorphisms
to finite $p$-groups, see \cite{Lyndon_Schupp}). Up to our knowledge
it is an open question if a product of several conjugacy classes
$[g_1]^F[g_2]^F...[g_k]^F$ is closed in the profinite topology in a
free group $F$. For pro-$p$ topology the answer is known: in the
2-generated free group $F=\langle x,y\rangle$ there exist elements
$g_1$ and $g_2$ such that $[g_1]^F[g_2]^F$ is not closed in any
pro-$p$ topology, even more, its closure is not contained in
$N(g_1,g_2)$, the normal subgroup generated by $g_1, g_2$, see
\cite{Howie84}. Precisely, in \cite{Howie84}  the following is
proven. If $g_1=x^{-2}y^{-3}$, $g_2=x^{-2}(xy)^5$ and $a=xy^2$, then
$a\not\in N(g_1, g_2)$. On the other hand, for any $i$, one has $a
\equiv w_i^{-1}g_1w_iv_i^{-1}g_2v_i \mod F_i$ for some $w_i,v_i\in
F$, where $F_0=F$ and $F_{i+1}=[F_i,F]$. So, for any homomorphism
$\phi$ from $F$ to a nilpotent group
$\phi(a)\in\phi([g_1]^F[g_2]^F)$ (it must be $\phi(F_i)=\{e\}$ for
some $i$ ). Since a finite $p$-group is nilpotent, $a$ belongs to
the closure of $[g_1]^F[g_2]^F$ in any pro-$p$ topologies, but
$a\not\in N(g_1, g_2)$. So, the same statement may be valid for the
profinite topology. Let $F$ be a free group and $X\subseteq F$.
 Let us denote the closure of $X$ in the profinite
topology on $F$ by $\overline{X}$.
\begin{conj}\label{conj2}
For a finitely generated free group $F$, there exists a sequence
$g_1,g_2,...,g_k \in F$ such that
$$
\overline{[g_1]^F[g_2]^F...[g_k]^F}\not \subseteq
N(g_1,g_2,....,g_k).
$$
\end{conj}
We hope, that some techniques of \cite{Liebeck,Nikolov1,Nikolov2}
could be useful for resolving the conjecture.
\section{Sofic groups}
\begin{defn}
Let $H$ be a finite group with a bi-invariant metric $d$. Let $G$ be
a group, $\Phi \subseteq G$ be a finite subset, $\epsilon > 0$, and
$\alpha >0$. A map $\phi:\Phi\to H$ is said to be a $(\Phi,
\epsilon, \alpha)$-homomorphism if:
\begin{enumerate}
\item For any two elements $a,b \in \Phi$, with $a \cdot b \in \Phi$,
  $d(\phi(a)\phi(b),\phi(a\cdot b))< \epsilon$
\item If $e_G \in \Phi$, then $\phi(e_G)=e_H$
\item For any $a \neq e_G$, $d(\phi(a),e_H)> \alpha$
\end{enumerate}
\end{defn}

\begin{defn} \label{def_sofic1}
The group $G$ is sofic if there exists $\alpha >0$ such that for any
finite set $\Phi \subseteq G$, for any $\epsilon > 0$ there exists a
$(\Phi,\epsilon, \alpha)$-homomorphism to a symmetric group $S_n$
with the  normalized Hamming metric $h$.
\end{defn}

Let us give an equivalent definition which is more widely used, see
\cite{es1,es2}.

\begin{defn} \label{def_sofic2}
The group $G$ is sofic if for any finite set $\Phi \subseteq G$, for
any $\epsilon > 0$ there exists a $(\Phi,\epsilon,
1-\epsilon)$-homomorphism to a symmetric group $S_n$ with the
normalized Hamming metric $h$.
\end{defn}

The following proposition, in fact, is contained in \cite{Pestov}.
\begin{prop}
The definitions~\ref{def_sofic1} and \ref{def_sofic2} are
equivalent.
\end{prop}
\begin{pf}
A bijection $[n^2]\leftrightarrow \{(i,j)\;:\;i,j\in [n]\}$
naturally defines an injective inclusion $S_n\times S_n\to S_{n^2}$.
(  $(i,j)f\times g=((i)f,(j)g)$  ). One can check that if $f,g\in
S_n$, then $(1-h_{n^2}(g\times g,f\times f))=(1-h_n(f,g))^2$. So, if
$\phi:G\to S_n$ is a $(\Phi,\epsilon_0,\alpha_0)$-homomorphism, then
$\phi\times\phi:G\to S_n\times S_n<S_{n^2}$ is a
$(\Phi,\epsilon_1,\alpha_1)$-homomorphism with
$\epsilon_1=2\epsilon-\epsilon^2$ and $\alpha_1=2\alpha-\alpha^2$.
Now repeating this operation one can make $\alpha_n$ as close to $1$
as one wants, then choose $\epsilon_0$ such that $\epsilon_n$ is as
small as one wants. (We suppose that $0<\epsilon<\alpha\leq 1$.)
\qed
\end{pf}

\section{w-sofic groups and profinite topology}

Definition~\ref{def_sofic1} of sofic groups appeals to the following
generalization:

\begin{defn}\label{def_w-sofic}
A group $G$ is called w-sofic if there exists $\alpha >0$ such that
for any finite set $\Phi\subset G$, for any $\epsilon>0$ there
exists a finite group $H$ with a bi-invariant metric $d$ and a
$(\Phi,\epsilon,\alpha)$-homomorphism to $(H,d)$.
\end{defn}

\begin{rem}
\begin{itemize}
\item In Definition~\ref{def_w-sofic} we do not ask the metric to be
normalized. So, $\alpha$ may be any fixed positive number.
\item It
is easy to see, that any sofic group is w-sofic.
\item The idea of
using bi-invariant metrics in the context of sofic groups appears in
\cite{Pestov}. We also discussed it with E. Gordon.
\end{itemize}
\end{rem}

\begin{thm}\label{th_i_1}
Let $F$ be a finitely generated free group and $N\vartriangleleft
F$. Then $F/N$ is w-sofic if and only if for any finite sequence
$g_1,g_2,...,g_k$ from $N$ one has
$\overline{[g_1]^F[g_2]^F...[g_k]^F}\subseteq N$.
\end{thm}
($\overline{X}$  denote the closure of $X$ in the profinite topology
on $F$.)
\begin{cor}
If there exists a non w-sofic group then there exists a finitely
presented non w-sofic group.
\end{cor}

\begin{pf}
It follows from the definition of w-sofic groups that if there
exists a non w-sofic group then there exists a finitely generated
non w-sofic group $G$. So, $G=F/N$ for a free group $F=\langle
x_1,x_2,...,x_n\rangle$. Then there exist $g_1,g_2,...,g_k\in N$
such that $\overline{[g_1]^F[g_2]^F...[g_k]^F}\not\subset N$. Now,
$N(g_1,g_2,...,g_k)\subseteq N$, where $N(g_1,g_2,...,g_k)$ is the
normal subgroup generated by $g_1,g_2,...g_k$. So,
$\overline{[g_1]^F[g_2]^F...[g_k]^F}\not\subset N(g_1,g_2,...,g_k)$
and the group $\langle x_1,x_2,...,x_n\;|\;g_1,g_2,...,g_k\rangle$
is not w-sofic. \qed
\end{pf}

\begin{conj}\label{conj1}
There exists a non w-sofic group.
\end{conj}

\begin{cor}
Conjecture~\ref{conj1} and Conjecture~\ref{conj2} are equivalent.
\end{cor}

To finish this section we present another point of view on
Theorem~\ref{th_i_1}. Let $F$ be a free group, $N\triangleleft F$
and $\tilde F$ be its profinite completion. It is clear that
$F<\tilde F$ and $N<\tilde F$, but in general $N \ntriangleleft
\tilde F$. Let $\hat N$ denote the minimal normal subgroup such that
$N<\hat N\triangleleft \tilde F$ and  $\tilde N$ denote the closure
of $N$ in $\tilde F$. It is easy to see that $\hat N\leq \tilde
N\triangleleft \tilde F$ and $\hat N=\tilde N$ iff $\hat N$ is
closed in $\tilde F$. The following corollary is a consequence of
Theorem~\ref{th_i_1} (with known facts about residually finite
groups).
\begin{cor}
The group $F/N$ is w-sofic iff $N=\hat N\cap F$ and the group $F/N$
is residually finite iff $N=\tilde N\cap F$.
\end{cor}

\begin{rem}
Let $S=[g_1]^{\tilde F} [g_1^{-1}]^{\tilde F} [g_2]^{\tilde F}
[g_2^{-1}]^{\tilde F} \cdots \- [g_k]^{\tilde  F} [g_k^{-1}]^{\tilde
F}$. Then $\hat N(g_1,g_2,...,g_k)=\bigcup_{n=1}^\infty S^n$. So,
$\hat N(g_1,g_2,...,g_k)$ is a closed set if and only if $\hat
N(g_1,g_2,...,g_k)=S^n$ for some $n$, due to the fact that $S$ is a
compact set, see \cite{Hartley}.
\end{rem}

\section{Bi-invariant metrics on finite groups}\label{sec_bi_metric}

Any bi-invariant metrics $d:G\times G\to \Rset$ on a group $G$ may
be defined as $d(a,b)=\|ab^{-1}\|$ where the ``norm"
$\|g\|=d(e_G,g)$ satisfies the following properties ($\forall g,h\in
G$):
\begin{enumerate}
\item $\|g\|\geq 0$ \label{p1}
\item $\|e_G\|=0$
\item $\|g^{-1}\|=\|g\|$ \label{p3}
\item $\|hgh^{-1}\|=\|g\|$ ($\|\cdot\|$ is a function of conjugacy classes) \label{p4}
\item $\|gh\|\leq \|g\|+\|h\|$ (it is, in fact, the triangle inequality for $d$.) \label{p5}
\end{enumerate}

\begin{rem}
In fact, properties 1-5 define only semimetric. However, $N=\{g\in
G\; :\; \|g\|=0\}$ is a normal subgroup. Since $\|\cdot\|$ is
constant on left (right) classes of $N$, it naturally defines metric
on $G/N$.
\end{rem}

Let us give examples of such metrics:
\begin{description}
\item[I.] The normalized Hamming metric on $S_n$ $\|x\|=\frac{|\{j\;|\;j\neq
jx\}|}{n}$.
\item[II.] One can use bi-invariant matrix metrics on a
representation of a group $G$. For example if one takes the
normalized trace norm (the normalized Hilbert-Schmidt norm) on an
unitary representations with character $\chi$ one gets:
$$
\|g\|_\chi=\sqrt{\frac{2\chi(e_G)-(\chi(g)+\chi^*(g))}{\chi(e_G)}}.
$$
(If $\chi$ is the fixed point character on $S_n$, then
$\|g\|_\chi=\sqrt{2h_n(e_{S_n},g)}$.)
\item[III.] The following construction will be used in the proof of
Theorem~\ref{th_i_1}. We will use ``conjugacy class graph" which is
an analogue of Cayley graph, where conjugacy classes are used
instead of group elements. Let $G$ be a (finite) group and $CG$ be
the set of all its conjugacy classes, let ${\cal C}\subset CG$. The
conjugacy graph $\Gamma(G,\C)$ is defined as follows:  its vertex
set $V=CG$ and for $x,y\in V$ there is an edge $(x,y)$ iff $x\subset
cy$ for some $c\in \C$. (The graph $\Gamma(G,\C)$ is considered as
an undirected graph.) See fig.1 for example.

\begin{figure}[h]
\center
\includegraphics[width=0.4\textwidth]{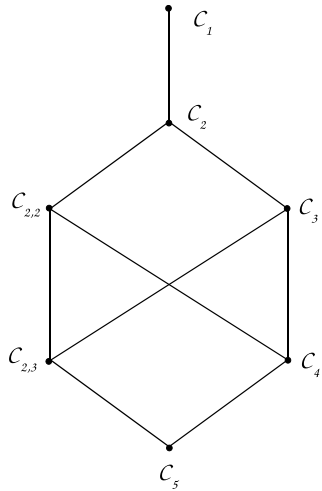}
\caption{The conjugacy graph $\Gamma(S_5,\{\C_2\})$, where the
vertices $\C_1, \C_2, \C_3, \C_4, \C_5, \C_{2,2},\C_{2,3}$ are the
conjugacy classes of $e_{S_5}, (1,2), (1,2,3), (1,2,3,4),
(1,2,3,4,5), (1,2)(3,4), (1,2)(3,4,5)$ respectively.}
\end{figure}

 Now
  define $\|g\|_\C$ to be
the distance from $\{e_G\}$ to $[g]^G$ in $\Gamma(G,\C)$, if
$[g]^G\in K(e_G)$, the connected component of $\{e_G\}$. For
$[g]^G\not\in K(e_G)$ let $\|g\|=\max\{\|x\|\;|\;x\in K(e_G)\}$. One
can check that $\|\cdot\|_\C$ satisfies \ref{p1}-\ref{p5} by
construction.
\end{description}

\section{Proof of Theorem~\ref{th_i_1}}

Theorem~\ref{th_i_1} is a direct consequence of the following
lemmata.

We will represent the elements of a free group $F$ by the reduced
words. For $w\in F$ let $|w|$ denotes the length of $w$.
\begin{defn}
Let $F$ be a free group, $N\triangleleft F$, $\delta>\epsilon>0$ and
$r\in\Nset$. Let $H$ be a finite group with a bi-invariant metric
$d$ and $\phi:F\to H$ be a homomorphism. We will say that $N$ is
$(r,\epsilon,\delta)$-separated by $\phi$ if for any $w\in F$,
$|w|\leq r$ we have the following alternative:
\begin{itemize}
\item $d(e_H,\phi(w))<\epsilon$ if $w\in N$,
\item $d(e_H,\phi(w))>\delta$ if $w\not\in N$.
\end{itemize}
In this case we call $N$ to be  $(r,\epsilon,\delta)$-separable.

We will say that a normal subgroup $N$ is finitely separable if
there exists $\delta>0$ such that for any $r\in\Nset$ any
$\epsilon>0$ the normal subgroup $N$ is
$(r,\epsilon,\delta)$-separable.
\end{defn}
\begin{lem}\label{lm_i_2}
Let $F$ be a finitely generated free group. Then $G=F/N$ is w-sofic
if and only if $N$ is finitely separable.
\end{lem}

\begin{pf}
Let us give two proofs of the lemma: the first one by using
nonstandard analysis and the second without it. A short introduction
to non standard analysis in a similar context could be found in
\cite{agg,Pestov}. Our use of nonstandard analysis is not essential:
one can easily rewrite the  proof on standard language (without
using ultrafilters), on the other hand the non-standard proof is
more algebraic and based on the following simple
\begin{claim}
Let $F$ be a free group and $H$ be a group, $N\triangleleft F$ and
$\E\triangleleft H$. Then the following are equivalent:
\begin{enumerate}
\item There exists a homomorphism $\phi:F\to H$ such that
$\phi(N)\subseteq \E$ and $\phi(F\setminus N)\cap \E=\emptyset$.
\item There exists an injective homomorphism $\tilde\phi:F/N\to
H/\E$.
\end{enumerate}
\end{claim}
\begin{pf}
Item 1 $\Rightarrow$ item 2. Let $\tilde\phi([w]_N)=[\phi(w)]_\E$.
It is  well-defined and injective by item~1.\\
Item 2 $\Rightarrow$ item 1. Let $F=\langle x_1,x_2,...,x_n\rangle$.
Let $[y_i]_\E=\tilde\phi([x_i]_N)$. Define $\phi:F\to H$ by setting
$\phi(x_i)=y_i$ ($F$ is free!). $\phi$ depends on the choice of
$y_i$, but for all $w\in F$ one has
$[\phi(w)]_\E=\tilde\phi([w]_N)$, (induction on  $|w|$). It proves
the claim. \qed
\end{pf}

{\bf Non-standard proof.}  Let $H$ be a hyperfinite group with an
internal bi-invariant metric $d$, then $\E=\{h\in
H\;|\;d(h,e)\approx 0\}$ is a normal subgroup (Since $d$ is
bi-invariant). Group $G=F/N$ is w-sofic iff there exists a
homomorphic injection $\tilde\phi:G\to H/\E$. On the other hand, $N$
is finitely separated iff there exists a homomorphism $\phi:F\to H$,
such that $\phi(N)\subseteq \E$ and $\phi(F\setminus
N)\cap\E=\emptyset$, this is equivalent to the existence of a
homomorphic injection $\tilde\phi:F/N\to H/\E$ by the claim. The
lemma follows.

{\bf Standard proof.}

$\Longrightarrow$ Let $G=F/N$ and $F=\langle
x_1,x_2,...,x_n\rangle$, let $\Phi=\{[w]_N\; | \; w\in F,\; |w|\leq
r\}$ and $\phi:\Phi\to H$ be a $(\Phi,\epsilon,\alpha)$-homomorphism
($[w]_N=wN$ denotes the right $N$-class of $w$). Define a
homomorphism $\tilde\phi:F\to H$ by setting
$\tilde\phi(x_i)=\phi([x_i]_N)$. It is enough to show that
$\tilde\phi$ will $(r,2r\epsilon,\alpha-2r\epsilon)$-separate $N$.
Using the inequality
$d(\phi([x_i]_N)^{-1},\phi([x_i^{-1}]_N))<\epsilon$ and induction on
$|w|$ one gets
$$
d(\tilde\phi(w),\phi([w]_N))< (2|w|-1)\epsilon.
$$
So, if $w\in N\cap \{w \in F \;:\; |w| \leq r\}$, then $[w]_N=e_G$
and
$$
d(\tilde\phi(w),e_H)< (2|w|-1)\epsilon < 2r\epsilon,
$$
if $w\not\in N$ then
$$
d(\tilde\phi(w),e_H)\geq
d(\phi([w]_N),e_H)-d(\tilde\phi(w),\phi([w]_N))>\alpha-2r\epsilon.
$$

$\Longleftarrow$ We have to construct a
$(\Phi,\epsilon,\alpha)$-homomorphism. Let
$\Phi=\{[w_1]_N,[w_2]_N,...,[w_m]_N\}$, such that $e_G$ (if exists)
is presented by the empty word, let $r=\max\{|w_1|,...,|w_m|\}$.
Choose $\phi$ to $(3r,\epsilon,\alpha)$-separate $N$, and define
$\tilde\phi:\Phi\to H$ as $\tilde\phi([w_i]_N)=\phi(w_i)$. (We
suppose, that all $[w_i]_N$ are different, so $\tilde\phi$ is
well-defined, although depends on the choice of $w_i$.) We claim
that $\tilde\phi$ is a $(\Phi,\epsilon,\alpha)$-homomorphism.
Indeed,
\begin{itemize}
\item $\tilde\phi(e_G)=e_H$,($e_G$ is represented by the empty word)
\item $d(\tilde\phi([w_i]_N,e_H)=d(\phi(w_i),e_H)>\alpha$
\item if $[w_i]_N[w_j]_N=[w_k]_N$ then
$$
d(\tilde\phi([w_i])\tilde\phi([w_j]),\tilde\phi([w_k]))=d(\phi(w_i)\phi(w_j),\phi(w_k))=
d(\phi(w_iw_jw_k^{-1}),e_H)<\epsilon,
$$
since $w_iw_jw_k^{-1}\in N $ and has length $\leq 3r$. \qed
\end{itemize}
\end{pf}

\begin{lem}
Let $F$ be a finitely generated free group and $N$ be its normal
subgroup. Then $N$ is finitely separable if and only if for any
$g_1,g_2,...,g_k\in N$ one has
$\overline{[g_1]^F[g_2]^F...[g_k]^F}\subseteq N$.
\end{lem}

\begin{pf}
$\Longrightarrow$ In order to prove that
$\overline{[g_1]^F[g_2]^F...[g_k]^F}\subseteq N$ it is enough to
prove that for any $w\not\in N$ there exists a homomorphism
$\phi:F\to H$ into a finite group $H$ such that $\phi(w)\not\in
[\phi(g_1)]^H[\phi(g_2)]^H...[\phi(g_k)]^H$.

Let $r>\max\{|w|,|g_1|,|g_2|,...,|g_k|\}$. Take $\phi:F\to H$ which
$(r,\alpha/k,\alpha)$-separates $N$, then $d(e_H,\phi(w))>\alpha$
and $d(e_H,\phi(g_i))<\alpha/k$. But if $\phi(w)\in
[\phi(g_1)]^H[\phi(g_2)]^H...[\phi(g_k)]^H$, then (by \ref{p4} and
\ref{p5} of Section~\ref{sec_bi_metric})
$$
d(e_H,\phi(w))\leq \sum_{i=1}^{k} d(e_H,\phi(g_i))<\alpha,
$$
a contradiction.

$\Longleftarrow$ We have to construct an
$(r,\epsilon,\alpha)$-separating homomorphism. Let $W_r=\{w\in
F,\;:\;|w|\leq r\}$ and integer $k>\alpha/\epsilon$. We can find a
homomorphism $\phi:F\to H$ into a finite group $H$, such that for
any $w\in W_r\setminus N$ and any $g_1,g_2,...,g_m\in W_r\cap N$,
$m\leq k$ one has $\phi(w)\not\in
[\phi(g_1)]^H[\phi(g_2)]^H...[\phi(g_m)]^H$. To define a metric on
$H$ let $\C=\{[\phi(g)]^H\;|\;g\in W_r\cap N\}$ and now let
$d(x,e)=\epsilon\|x\|_\C$, where $\|\cdot\|_\C$ is defined in item
III of section~\ref{sec_bi_metric}. Now one can check that
$\phi:F\to H$ will $(r,\epsilon',\alpha)$-separate $N$, for any
$\epsilon'
> \epsilon$.\qed
\end{pf}



\end{document}